# Perfect cuboid, primitive Pythagorean triples and Eulerian parallelepipeds. Dynamics of construction

*Aleshkevich Natalia V.*

*Peter the Great St. Petersburg Polytechnic University*

**Abstract**

One unsolved mathematical problem remains the perfect cuboid problem. A perfect cuboid is a rectangular parallelepiped whose edges, face diagonals and space diagonal are all expressed as integers. No such cuboid has yet been discovered and its existence has also not been proven. This paper shows a proof of the non-existence of a perfect cuboid.

# Introduction

Before starting to work on proving the existence or non-existence of a perfect cuboid, we will consider such objects as primitive Pythagorean triples, general Pythagorean triples, and Eulerian parallelepipeds.

One of the objects of our consideration are the Pythagorean numbers, also called Pythagorean triples – triples $(x, y, a)$ of natural numbers satisfying the Pythagorean equation

$$x^2 + y^2 = a^2.$$

The Pythagorean theorem is a fundamental geometric statement: in any right triangle, the area of a square built on the hypotenuse is equal to the sum of the areas of squares built on the legs.

The general solution are the following formulas [1]:

$$y = 2mn; x = m^2 - n^2;  \quad a = m^2 + n^2.$$

These formulas describe exactly once every Pythagorean triple $(x, y, a)$, satisfying the condition $\text{GCD}(x, y, a) = 1$. This means that all sides of the Pythagorean triangle are expressed by relatively prime numbers. This triple of numbers is called a primitive Pythagorean triple.

In any primitive Pythagorean triple one of the legs is an even number and the other is an odd number. In this case the hypotenuse $a$ is an odd number. Without loss of generality, we will assume that $x$ is odd and $y$ — even. Under these constraints we can get all primitive Pythagorean triples and only them.

The parameters $m$ and $n$, forming primitive Pythagorean triples, were obtained from very abstract considerations and are not related to each other; that is, they are independent.

The task was to find a geometric interpretation of generation of primitive Pythagorean triples; and to, based on the received interpretation, determine the order on the set of primitive Pythagorean triples, their properties, and quantitative estimates.

The next object of our consideration are Eulerian parallelepipeds [2] with integer edges and face diagonals. They are described by the following system of equations:

$$\begin{cases} x^2 + y^2 = a^2 \\ y^2 + z^2 = b^2 \\ x^2 + z^2 = c^2 \end{cases} \quad (1)$$

The minimal parallelepiped with edges $(117, 44, 240)$ was found in 1719 by the German mathematician Paul Halcke [3]. Leonhard Euler proposed a particular solution for finding edges and face diagonals of a parallelepiped. Therefore, parallelepipeds were called Eulerian.

No general algorithm for constructing such parallelepipeds had ever been proposed till today. This algorithm is presented in this paper.

And the last object of our study is a perfect cuboid: it is an Eulerian parallelepiped with an integer main (space) diagonal $d$. In this case, to the system of equations (1) the next equation is added:

$$x^2 + y^2 + z^2 = d^2. \quad (2)$$

The paper proves the impossibility of constructing a perfect cuboid.

**Dynamics of construction of primitive Pythagorean triples**

We will consider the construction of consecutive squares, starting with the square of one. The formula of such a construction:

$$(n + 1)^2 = n^2 + 2n + 1.$$

To the constructed square with side $n$, we can add a figure whose area is equal to twice the value of the side of the square plus 1: $2n + 1$. This figure, called a gnomon, builds the original square to a larger square; the side of which will be equal to $2n + 1$. The thickness of the gnomon will be equal to 1. By constructing $n$ such consecutive gnomons, we can construct a square with side $n + k$. We can combine consecutive gnomons with a thickness equal to 1 into one common gnomon with a thickness equal to $k$.

$$x^2 + G = a^2.$$

We need to build a gnomon that is equal in area to some square: $G_y = y^2$

Then we come to the equation

$$x^2 + y^2 = a^2.$$

**Building squares and their sum using a generating square**

We show the construction of squares of a primitive Pythagorean triple using a generating square with side $S = 2tl$ [4]. Here $\gcd(t, l) = 1$. We assume that $l$ is odd and $t$ is of any parity.

Without loss of generality, we will build a square with an even side.

We increase the side $S$ of the generating square by $2t^2$ and build a larger square with side $y = S + 2t^2$ (Fig. 1). In this case, we obtain a gnomon $G$ with a thickness $2t^2$, placed on the generating square.

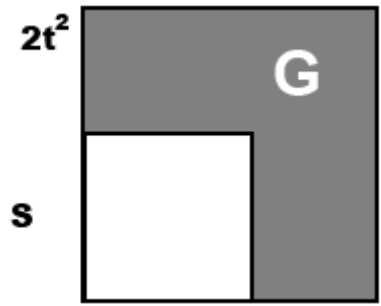

Figure 1

Next, we increase the $y$ side by the value $l^2$. Concurrently we extend the side of the gnomon by the same value $l^2$. At both ends of the gnomon we will have identical rectangles with an area $2t^2 \times l^2$ (Fig. 2).

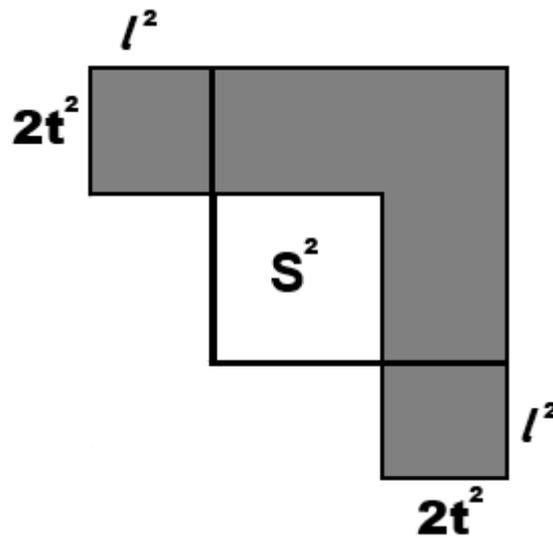

Figure 2

The total area of both rectangles is equal to the area of the generating square with side $S$. If we remove the generating square we obtain the gnomon $G_y$; which is equal in area to the square with side $y$ (Fig. 3).

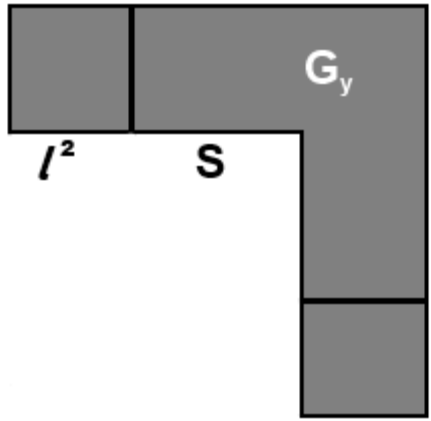

Figure 3

Gnomon $G_y$ is placed on a square with a side:

$$x = S + l^2 = 2tl + l^2 = l(l + 2l) \text{ (Fig. 4).}$$

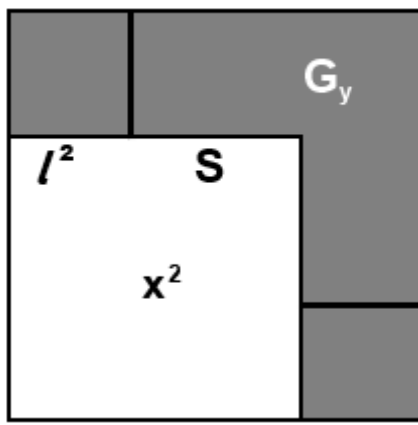

Figure 4

The outer side of the gnomon $G_y$ is equal to the hypotenuse:

$$a = S + 2t^2 + l^2 = (l + t)^2 + t^2.$$

The sum of two squares can be represented as one of the squares and a gnomon placed on it, which is equal in area to the second square. This representation is symmetrical (Fig. 5, 6):

$$x^2 + G_y = y^2 + G_x.$$

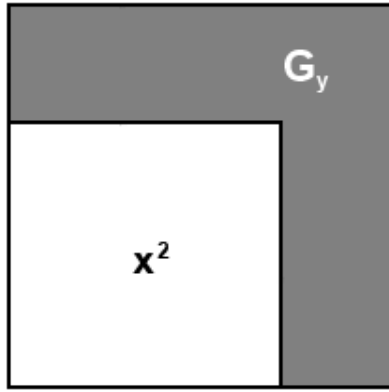

Figure 5

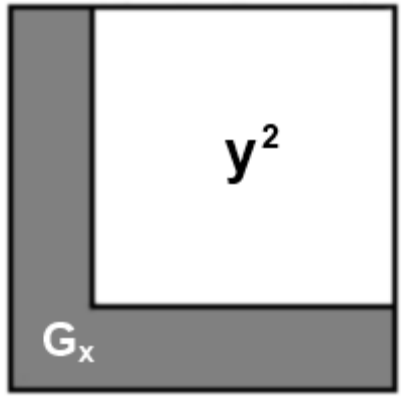

Figure 6

The total square in the form of both gnomons is shown in Fig. 7. Here, the larger gnomon absorbs the smaller gnomon. We will call this representation of gnomons connected gnomons. Both gnomons have a common outer side equals to the hypotenuse $a$. Thus, we have the following relation:

$$x^2 + G_y = y^2 + G_x = G_x + G_y = a^2.$$

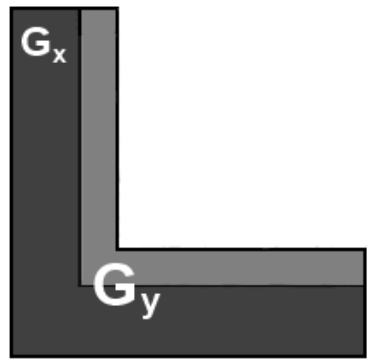

Figure 7

## Expression of the parameters $m$ and $n$ through the partitioning parameters of the side of the generating square

Expression of the parameters $m$ and $n$ through the partitioning parameters of the side of the generating square was considered earlier. [4]

The side of the even square is represented as $y = 2t(l + t)$.

We decided on $l$ is an odd number. It can be either 1 or the product of odd prime cofactors in the corresponding powers. The number $t$ can be any parity. It can be either 1 or the product of prime cofactors in the corresponding powers. We obtain that the factors $t$ and $(l + t)$ have different parity. And, indeed, if $t$ is even, then the contents of the bracket are odd and vice versa.

The transition from our notation to the generally accepted in terms of $m$ and $n$:

$$m = l + t; n = t.$$

In this case, the even square will have a side $y = 2mn$. The total square will have a side $m^2 + n^2$. In fact,

$$a = 2tl + 2t^2 + l^2 = t^2 + 2tl + l^2 + t^2 = (t + l)^2 + t^2.$$

The square with an odd side (odd square) will have a side $m^2 - n^2$. In fact,

$$x = 2tl + l^2 = 2tl + l^2 + t^2 - t^2 = (t + l)^2 - t^2.$$

As a result of our investigation were revealed the algebraic meaning that underlie the choice of the parameters $m$ and $n$:

$m$ - is the sum of the values of two subsets of the partition;

$n$ - is the value of one subset of the partition, it necessarily includes the factor $2^{\alpha_0 - 1}$ with $\alpha_0 > 1$.

Consequently, the parameters $t$ and $l$ also definite the primitive Pythagorean triple unambiguously, as well as the parameters $m$ and $n$.

**Setting the order on a set of primitive Pythagorean triples**

Setting the order on a set of primitive Pythagorean triples was considered earlier. [5]

Our construction is based on an even square (a square with an even side) generating the sum of two squares, which is also a square. The side of the generating square $S$ is represented as a product of $S = 2tl$.

We can consider a sequence of generating squares, starting with a square with side $S = 2$, and then move with a common difference of $2$.

Represent the side of the generating square $S$ as: $S = 2 \times 2^{\alpha_0 - 1} p_1^{\alpha_1} p_2^{\alpha_2} \ldots p_r^{\alpha_r}$, where $p_i$ are prime odd cofactors.

The amount of $L(S)$ partitions of the product into two groups of factors $t$ and $l$ depends on the amount of $r$ odd prime cofactors without taking into account their powers and is equal to $L(S) = 2^r$; that is, equal to the sum of the binomial coefficients for the row with the number $r$ in Pascal's triangle. Thus, for each $S$ there are $2^r$ primitive Pythagorean triples.

The entire set of primitive Pythagorean triples can be constructed according to the sequential growth of the known parameters. [5]

These parameters are the side of the generating square $S$ and factors $t$ and $l$. Since the parameters $t$ and $l$ are mutually related, we will choose one of them for ordering, namely $t$. Thus, the order is set by two parameters. One parameter external is the side of the generating square. Side $S$ is even number. $S$ starts with $2$ and goes in increments of $2$. The internal parameter $t$ is the partition element of the side of the generating square $S = 2tl$. The element $t$ starts with the minimum value corresponding to the parameter $S$, and then increases to the maximum value within $S$. The associated element $l$, starting from the maximum, decreases at the same time. Both elements are formed from the cofactors of the number $S$. The number $t$ can be of any parity. The number $l$ is odd. Wherein

$$GCD(t, l) = 1.$$

Formulas for obtaining elements of a primitive Pythagorean triple:

$$y = S + 2t^2;$$

$$x = S + l^2;$$

$$a = S + 2t^2 + l^2.$$

According to the parameters $S, t(S)$ ordered tables of primitive Pythagorean triples can be constructed. In this case, the $N$ ordinal number of the first level is equal to

$$N = S/2.$$

The sequence number $n$ of the second level changes from 1 to $L(S)$ within $S$, where

$$L(S) = \sum_{i=0}^{r} C_r^i = 2^r,$$

r - is the number of prime odd cofactors without taking into account their powers included in the product for $S$.

The table is presented in Appendix 1. Where is a fragment of the beginning of the set for the values $S = 2 \div 100$. Accordingly $N = 1, 2, ..., 50$.

Setting the order allows to build different algorithms when using primitive Pythagorean triples.

## Primitive Pythagorean triples and their representation via arithmetic progressions

Primitive Pythagorean triples and their representation via arithmetic progressions also was considered earlier. [5]

Let's imagine a primitive Pythagorean triple in the form of a square and a gnomon placed on it. Take a square with an odd side $x$. Then the area of the gnomon can be represented as the sum of an arithmetic progression with the first term $2x + 1$. Each subsequent member will be two units larger than the previous one. The number of such terms in the arithmetic progression is equal to the thickness of the gnomon

$$T_y = 2t^2.$$

Now take a square with an even side $y$. Then the area of the gnomon built on it can be represented as the sum of an arithmetic progression with the first term equal to $2y + 1$. Each subsequent member will be two units larger than the previous one. The number of such terms in the arithmetic progression is equal to the thickness of the gnomon

$$T_x = l^2.$$

We take two connected gnomons. For $y < x$, all the terms of the gnomon $G_y$, and their number is $T_y$, will be equal, respectively, to the last terms in the arithmetic progression representing $G_x$. And, conversely, for $x < y$, all the terms of the gnomon $G_x$, and their number is $T_x$, will be equal, respectively, to the last terms in the arithmetic progression representing $G_y$.

This representation in the form of an arithmetic progression of each gnomon fully corresponds to the picture of the gnomon absorbing a larger area of the connected gnomon of a smaller area.

The sum of the terms of the arithmetic progression is equal to the square of the corresponding leg.

The sum of the terms of both arithmetic progressions is equal to the square of the hypotenuse.

The middle term $s_x$ of the arithmetic progression (arithmetic mean) describing the gnomon $G_x$ is equal to the sum of the arithmetic progression divided by the number of its terms:

$$s_x = \frac{x^2}{T_x} = \frac{l^2(l+2t)^2}{l^2} = (l+2t)^2.$$

The middle term $s_y$ of the arithmetic progression (arithmetic mean) describing the gnomon $G_y$ is equal to the sum of the arithmetic progression divided by the number of its terms:

$$s_y = \frac{y^2}{T_y} = \frac{4t^2(l+t)^2}{2t^2} = 2(l+t)^2.$$

The formula for the first term of the arithmetic progression:

$$s1 = s - T + 1.$$

Since $s1_y = 2x + 1$ and $s1_x = 2y + 1$, therefore we can find the values for $x$ and $y$:

$$x = \frac{s1_y - 1}{2} = \frac{s_y - T_y}{2} = \frac{2(l+t)^2 - 2t^2}{2} = (l+t)^2 - t^2 = l^2 + 2tl$$
$$= l(l+2t).$$

$$y = \frac{s1_x - 1}{2} = \frac{s_x - T_x}{2} = \frac{(l+2t)^2 - l^2}{2} = \frac{4t^2 + 4tl}{2} = 2t(l+t).$$

It follows that each gnomon, equal in area to the square of one of the legs, is placed on the square of the second leg from the primitive Pythagorean triple. Thus, there is a mutual mapping of the legs:

$$x \to y;\ y \to x;$$

$$x = l(l + 2t) \leftrightarrow y = 2t(l + t).$$

For mapping $x \to y$, we use the formula

$$\frac{(l + 2t)^2 - l^2}{2}.$$

For mapping $y \to x$, we use the formula

$$(l + t)^2 - t^2.$$

Connected gnomons match with the last terms of arithmetic progressions. The last term of the arithmetic progressions $s_n$ is equal to the sum of the middle term of the arithmetic progression and the corresponding amount of terms in this progression minus one:

$$s_n = s_x + T_x - 1 = s_y + T_y - 1.$$

In this case, the last term is equal to $s_n = 2a - 1$. Hence the equality for the hypotenuse $a$ follows:

$$a = \frac{s_n + 1}{2};$$

$$a = \frac{s_x + T_x}{2} = \frac{(l + 2t)^2 + l^2}{2} = \frac{2l^2 + 4t^2 + 4lt}{2} = 2lt + 2t^2 + l^2;$$

$$a = \frac{s_y + T_y}{2} = \frac{2(l + t)^2 + 2t^2}{2} = \frac{2l^2 + 4t^2 + 4lt}{2} = 2lt + 2t^2 + l^2.$$

We substitute $S = 2lt$ into the equations and have:

$$a = S + 2t^2 + l^2 = x + 2t^2 = y + l^2.$$

The last equation corresponds to the value of the hypotenuse $a$, constructed by means of a generating square with side $S$.

Thus, using the concept of arithmetic progression to describe the connected gnomons of a primitive Pythagorean triple, we obtained that the connected gnomons $G_x$ and $G_y$ uniquely represent the primitive Pythagorean triple $(y, x, a)$.

**Transformation of Gnomons**

A gnomon can be transformed into another gnomon while preserving its area.

The transformation of gnomons leads to a restructuring of their structure: their thickness and the value of the middle term of the arithmetic progression describing the gnomon, change in a coordinated manner so that the size of the gnomon area is conserved. As the thickness of the gnomon decreases, the middle term of the arithmetic progression increases, and vice versa, as the thickness of the gnomon increases, the middle term of the arithmetic progression decreases. During the transformation the gnomon, equal in area to the square of the first leg, is placed on another square that is different from the square of the second leg from the initial primitive Pythagorean triple. In this case, the side of the total square also changes.

The amount of transformations of the gnomon depends on the prime factorization of the square's side that the gnomon represents. During the transformation, the gnomon can form some amount of new primitive Pythagorean triples. This applies to each connected gnomon. Specifically, this pair of legs from the primitive Pythagorean triple disintegrates when at least one gnomon is transformed.

The amount of new pairs of legs is determined by the following ratio. For each leg, it depends on the number of its cofactors. For an even leg, it depends on the possible number of partitions of its cofactors into groups $2t$ and $(l + t)$. For an odd leg, it depends on the possible number of partitions of its cofactors into groups $l$ and $(l + 2t)$. At the same time, the parameter $S$ changes for the transformed gnomon. Specifically, the gnomon after transformation forms another primitive Pythagorean triple. For the leg $y$, other legs are selected and accordingly, other hypotenuses are obtained; the same applies for the leg $x$. As many transformations of the gnomon are possible, so many new primitive Pythagorean triples can be built. The parameters $l$ and $t$ inside each pair of legs are relatively prime. Thus, we can determine the amount of identical legs $y$ included in different rows of the table of primitive Pythagorean triples. The same applies for the legs $x$.

The amount of different $S$ for the same leg

$$y = 2t(l + t)$$

is determined by the amount of different $t$ that can be substituted into the formula for $y$. And this amount is determined by the sum of two Stirling numbers of the second kind $S(n, k)$. Where $n$ is the number of prime factors of $y$, without taking in account their powers, and $k = 1, 2$:

$$S(n, 1) + S(n, 2) = 1 + 2^{n-1} - 1 = 2^{n-1} \quad (3)$$

The number of different $S$ for the same leg $x = l(l + 2t)$ is determined by the number of different $l$ that can be substituted into this formula for $x$. And this amount is determined by the sum of two Stirling numbers of the second kind $S(n, k)$. Where $n$ is the number of prime factors of $x$ without powers and $k = 1, 2$. The formula is identical to (3).

By definition [6], the Stirling numbers of the second kind $S(n, k)$ are the number of ways to partition a set of $n$ objects into $k$ non-empty subsets, if $n = 1, 2 \ldots; k = 1, 2 \ldots, n$.

In our case, $n$ is the number of prime factors, without taking into account their powers, in the product for each particular leg ($y$ or $x$). According to the formulas $x = l(l + 2t)$ and $y = 2t(l + t)$, these products can be divided into one and two parts, that is, $k = 1, 2$.

**General Pythagorean triples**

We considered general Pythagorean triples earlier. [5]

Multiply all the elements of a primitive Pythagorean triple by an integer coefficient $k$.

Let's construct a primitive Pythagorean triple $(x, y, a)$. Let's draw it as a square with side $a$. Inside this square, in the upper right row, we place a square with the side $y$. Let's add gnomon $G_x$ to the inner square. The area of gnomon $G_x$ is equal to the area of a square with side $x$. The thickness of the gnomon is $T_x = l^2$. We place the square $a^2$ in the cells of the square lattice with the side $ka$ (Fig. 8).

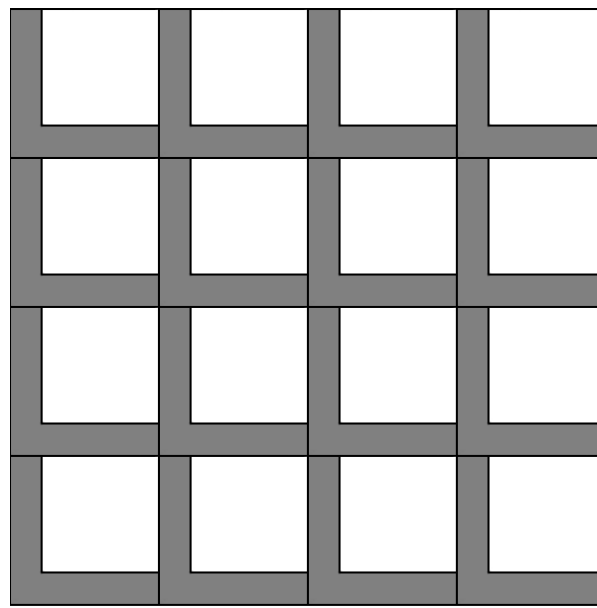
Figure 8

Let's put together all the squares with the $y$ side on the right, and on the left and at the bottom we will draw the total gnomon from the gnomons of each square lattice cell (Fig. 9). We see that the thickness of the gnomon $G_{kx}$ has become equal to $T_{kx} = kl^2$. The area of the summing square $a^2$ has increased $k^2$ times. The area of the square with the side $y$ has increased $k^2$ times. Consequently, the area of the gnomon also increased $k^2$ times. Thus, we have constructed a general Pythagorean triple $(kx, ky, ka)$.

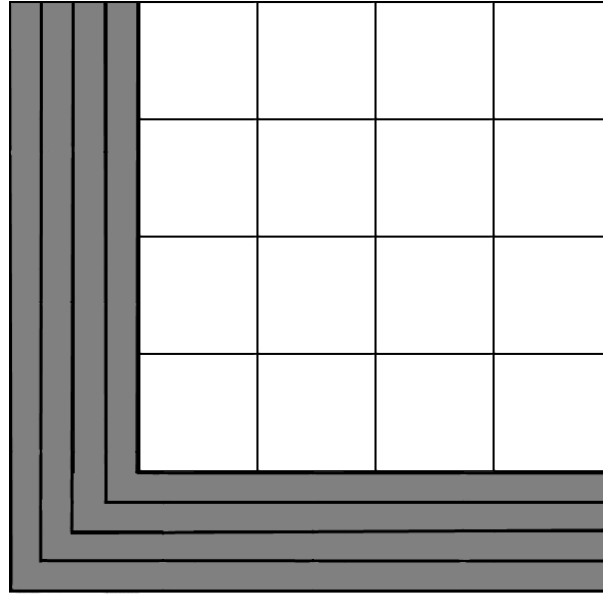

Figure 9

The construction will be similar, if in the total square with side $a$, we place a square with side $x$ and place a gnomon $G_y$ on it. The thickness of the gnomon is $T_y = 2t^2$. After assembling the squares and gnomons separately on a square lattice, we will construct a general Pythagorean triple $(kx, ky, ka)$. The area of the square with side $x$ will increase by $k^2$ times. The thickness of the gnomon $G_{ky}$ will increase by $k$ times and become equal to $T_{ky} = 2kt^2$.

The lengths of the outer sides of both gnomons will be equal to $ka$.

Thus, when we multiply all the elements of the Pythagorean triple by an integer coefficient $k$, the thickness of each gnomon in the corresponding constructions increases by $k$ times.

## Algorithm for constructing Eulerian parallelepipeds

It is easy to see that from each Pythagorean triangle we can obtain a rectangle whose sides and diagonals are expressed in natural numbers; and vice versa, any rectangle of this kind generates a Pythagorean triple. A rectangular parallelepiped contains three such original rectangles in the planes $XY, YZ, XZ$. That is, we have to solve a system of three equations (formula 1).

In formula 1, $x, y, z$ are the edges of the parallelepiped and $a, b, c$ are its face diagonals.

We will work with a table of primitive Pythagorean triples (appendix 1). We will sequentially iterate over the parameters $S$, starting with $S = 2$. And inside the block $S$, we will iterate over the parameter $t(S)$. A specific pair of parameters $S, t(S)$ uniquely defines a primitive Pythagorean triple $(x, y, a)$.

It follows from the first equation of the system (formula 1) that the numbers $x$ and $y$ have different parity, since they are the legs of a primitive Pythagorean triple. Without loss of generality, we assume that $x$ is an odd number and $y$ is an even number. Then in the third equation of the system (formula 1) we have $z$ - an even number; since both legs cannot be odd. It follows that the second equation contains two even legs, that is, it builds a general Pythagorean triple. Therefore $y$ and $z$ are scaled legs, meaning they have a common multiplier. At the same time the $x$ and $z$ legs can have common odd multipliers, that is, they can also build a general Pythagorean triple, or they can be relatively prime.

For each pair of legs from a primitive Pythagorean triple we can consider the possibility of constructing a third leg to build an Eulerian parallelepiped.

In general, the possibility of finding the third edge of a parallelepiped, or a common leg for a pair of legs from a primitive Pythagorean triple, depends on the representation of the squares of legs of a primitive Pythagorean triple in the form of gnomons transformed in such a way that each gnomon is placed on the same square. The side of this square will be the desired third edge of the parallelepiped (Fig. 10).

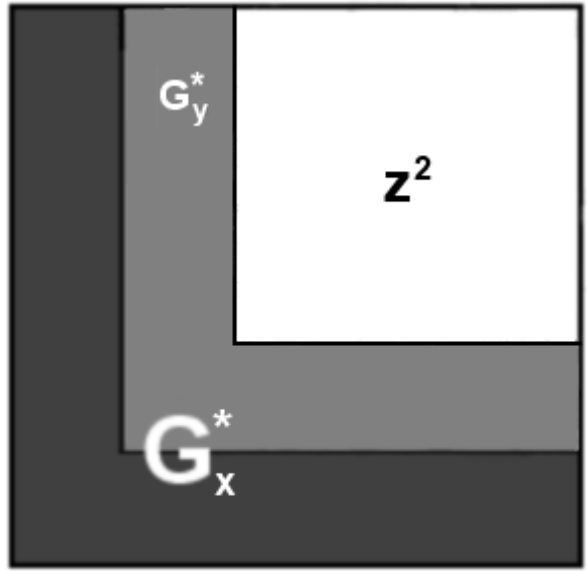

Figure 10

Two parameters change simultaneously in a gnomon described by an arithmetic progression; during its transformation: the number of terms of the arithmetic progression (or the thickness of the gnomon $T$) and the middle term of the arithmetic progression $s$. The area of the gnomon remains unchanged.

After transformation the middle terms of arithmetic progressions are equal, respectively:

$$s_y^* = \frac{y^2}{2g_y t^2}; \quad s_x^* = \frac{x^2}{g_x l^2}.$$

The thickness of the gnomons (or the number of terms of the arithmetic progression) are equal to:

$$T_y^* = 2g_y t^2; \quad T_x^* = g_x l^2,$$

where $g_x$ and $g_y$ are the transformation coefficients of the gnomons $G_x$ and $G_y$, respectively.

The first terms of the arithmetic progression $s1$ in both transformed gnomons $G_x^*$ $G_y^*$ must be equal to the same number:

$$s1_x = s1_y = 2z + 1.$$

The side of the square $z$ in this case will be equal for the gnomon $G_y^*$:

$$z = \left(\frac{y^2}{2g_y t^2} - 2g_y t^2\right)/2;$$

and for the gnomon $G_x^*$:

$$z = \frac{\frac{x^2}{g_x l^2} - g_x l^2}{2}.$$

Under this condition, both gnomons will be placed on the same square. The side $z$ of this square will be the desired third edge of the parallelepiped.

To find a new leg $z$ for an even leg $y$, since $z$ is also an even leg, we will use a square lattice.

We describe the general constructing of the side of the square lattice for the case when the legs $x$ and $z$ also have a common multiplier which is odd.

First we select a primitive Pythagorean triple from the table of primitive triples:

$$(S, t, l); \quad (x, y, a);$$

$$S = 2tl; \quad y = 2tl + 2t^2 = 2t(l + t); x = 2tl + l^2 = l(2t + l).$$

$$a = 2tl + 2t^2 + l^2 = (t + l)^2 + t^2.$$

Then we work with the legs $x$ and $y$. In a parallelepiped, the legs are edges. The hypotenuse $a$ is a face diagonal constructed on the edges $x, y$.

We need to find a common edge $z$ for edges $x$ and $y$, such that the following equalities hold together:

$$y^2 + z^2 = b^2;$$

$$x^2 + z^2 = c^2.$$

We start the search for the third leg, also known in this representation as the third edge for the parallelepiped, by working with an even leg.

Imagine an even leg $y$ as a product $y = 2t(l + t)$.

For an even $t$, we have an odd multiplier in brackets. For an odd $t$ we have an even multiplier in brackets. In this case, $y$ will have the form $y = 4i, i \in \mathbb{N}$.

Represent $y$ as a product of $y = k_1 m_1$, where $k_1$ is a coefficient, and $m_1$ is a leg that is part of one or more primitive Pythagorean triples depending on the number of its various prime factors. The coefficient $k_1$ will always be an even number of the form $k_1 = 4n,$ where $n \in \mathbb{N}$ (see below). For each selected coefficient $k_1$, we have the corresponding truncated leg $m_1$. When choosing a coefficient the truncated leg cannot be equal to 1 since there is no such leg for primitive Pythagorean triples.

The leg $m_1$ can be either an even or an odd number.

In general, one leg can be included in several different primitive Pythagorean triples, but the list of these triples is limited by the number of possible $S_i$; since $S_i$ is always smaller than a leg in numerical value. That is, the list of possible $S_i$ has a limit on the size of the list. For an even leg, the larger is $t$ and the smaller is $S$. The value $S$ must have $t$ as its multiplier. For an odd leg, the larger is $l$ and the smaller is $S$. The value $S$ must have $l$ as its multiplier:

$$y = S + 2t^2;$$

$$x = S + l^2.$$

Consider the first case where the leg $m_1$ is an odd number. We make a list of possible $l_i$. If $m_1$ is a prime number, then $l = 1$. If $m_1$ is a composite number, then we represent it as a product of prime factors in the corresponding powers:

$$m_1 = r_1^{\alpha_1} r_2^{\alpha_2} \ldots r_k^{\alpha_k}.$$

In this case, parameter $l_i$ can be equal to 1 or any combination of factors included in $m_1$, on condition that

$$l_i < \frac{m_1}{l_i}.$$

For each $l_i$ on the list, we find the leg $m_3$ from the corresponding primitive Pythagorean triple for the leg $m_1$. The leg $m_3$ can be calculated by the formula:

$$m_3 = \frac{(\frac{m_1}{l_i})^2 - l_i^2}{2}.$$

If the leg $m_1$ is an even number, then we make a list of possible $t_i$. In this case, the parameter $t_i$ can be equal to 1 or any combination of factors included in $m_1$, on condition that

$$t_i < \frac{m_1}{t_i}.$$

For each $t_i$ on the list, we find the leg $m_3$ from the corresponding primitive Pythagorean triple for the leg $m_1$. The leg $m_3$ can be calculated by the formula:

$$m_3 = (\frac{m_1}{t_i})^2 - t_i^2.$$

Under a fixed coefficient of truncated $k_1$ the number of different pairs of legs $(m_1, m_3)$ is determined by the formula (formula 3).

The leg $m_3$ will be the side of the inner square in the square lattice with a side $z$. The square of the leg $m_1$ is placed on this inner square in the form of a gnomon. Multiplying the truncated legs by the coefficient $k_1$, we obtain

$$y = k_1 m_1;$$

$$z = k_1 m_3.$$

The gnomon, representing the square of the leg $m_1$, has a thickness $T_{m_1} = l_1^2$, if $m_1$ is an odd number or it has a thickness $T_{m_1} = 2t_1^2$, if $m_1$ is an even number. After assembling the gnomon $G_y^*$ on a square lattice, its thickness will be equal to either $T_y^* = k_1 l_1^2$, or $T_y^* = 2k_1 t_1^2$.

As a result of the construction, we found a complete list of possible candidates $z_i$ to be the third edge of the parallelepiped.

Thus we can construct the following equations for the selected parameters $k_1, l_1, t_1$:

$$y^2 + z_i^2 = b_i^2,$$

where $i$ - is the number of the element in the lists for $l_1$ and $t_1$.

We make a similar action for leg $x$.

In general, the odd leg $x$ is represented by the formula:

$x = l(l + 2t)$.

Represent $x$ as a product of $x = k_2 m_2$, where $k_2$ is a coefficient, and $m_2$ is a leg that is part of one or several primitive Pythagorean triples (depending on the number of its various prime factors). The coefficient $k_2$ is an odd number. For each selected coefficient $k_2$, we have the corresponding truncated odd leg $m_2$.

When choosing a coefficient the truncated leg cannot be equal to 1, since there is no such leg for primitive Pythagorean triples.

We make a list of possible $l_i$. If $m_2$ is a prime number, then $l = 1$. If $m_2$ is a composite number, then we represent it as a product of prime factors in the corresponding powers:

$$m_2 = r_1^{\alpha_1} r_2^{\alpha_2} \ldots r_k^{\alpha_k}.$$

In this case, the parameter $l_i$ can be equal to 1 or any combination of factors included in $m_1$, under the condition

$$l_i < \frac{m_2}{l_i}.$$

For each $l_i$ on the list, we find the second leg $m_4$ from the corresponding primitive Pythagorean triple for the leg $m_2$. The leg $m_4$ can be calculated by the formula:

$$m_4 = \frac{(\frac{m_2}{l_i})^2 - l_i^2}{2}.$$

The number of different pairs of legs $(m_2, m_4)$ under a fixed coefficient of truncation $k_2$, depends on the number of prime factors without taking into account the powers of leg $m_2$ and is determined by the formula (formula 3).

The leg $m_4$ is the side of the inner square in the square lattice with the side $z$. The square of the leg $m_2$ is placed on this inner square in the form of a gnomon. When multiplying the truncated legs by the coefficient $k_2$, we obtain:

$$x = k_2 m_2;$$

$$z = k_2 m_4.$$

The gnomon, representing the square of the leg $m_2$, has a thickness $T_{m_2} = l_2^2$. After assembling the gnomon $G_x^*$ on a square lattice, we have $T_x^* = k_2 l_2^2$.

As a result of the construction, we found a complete list of possible candidates for $z_i$ to be the third edge of the parallelepiped.

Thus we can construct the following equations for the selected parameters $k_2, l_2,$:

$$x^2 + z_i^2 = c_i^2,$$

where $i$ - is the number of the element in the list for $l_2$.

Therefore, for each leg from the primitive Pythagorean triple $(x, y, a)$, we find a list of applicants for leg $z$. To build an Eulerian parallelepiped for the legs $(x, y)$, we need two values of $z_i$ from different lists to match. If there is no match then it is impossible to construct an Eulerian parallelepiped for this particular pair from the table of primitive Pythagorean triples.

We will find formal conditions for matching of two elements from different lists of applicants. The legs $x$ and $y$ are relatively prime by definition. The truncation coefficients of the legs $k_1$ and $k_2$ are also relatively prime. The truncated legs $m_1$ and $m_2$ are also relatively prime. Consider the legs $m_3$ and $m_4$ paired to them.

Imagine $m_3$ as a product of two groups of factors:

$$m_3 = k_2 q.$$

In this case, for an even leg $y$, we have $z = k_1 m_3 = k_1 k_2 q$.

Imagine $m_4$ as a product of two groups of factors:

$$m_4 = k_1 q.$$

In this case, for an odd leg $x$, we have $z = k_2 m_4 = k_2 k_1 q$.

When these conditions are met, we assert that coefficient of truncation $k_1$ will always be an even number of the form $k_1 = 4n$, where $n \in \mathbb{N}$. Otherwise, the truncated leg $m_1$ is always even, and the paired with it leg $m_3 = k_2 q$ is odd. And hence $q$ is also odd, and hence $z = k_1 k_2 q$ is also odd. This contradicts the condition from formula 1 (under this condition we construct only an even $z$).

If these conditions are met, we will obtain a common leg $z$ for legs $x$ and $y$ from a row in the table of primitive Pythagorean triples. Thus, we have found the necessary and sufficient conditions for the construction of the third leg:

$$\text{GCD}(m_3, x) = k_2;$$

$$\text{GCD}(m_4, y) = k_1;$$

$$\frac{m_3}{k_2} = \frac{m_4}{k_1} = q.$$

Under these conditions, we can simplify the algorithm for constructing the leg $z$ as follows.

In our construction $\text{GCD}(m_3, x) = k_2$.

$$m_2 = \frac{x}{k_2}.$$

Symmetrically for the case when the construction begins with leg $x$:

$$\text{GCD}(m_4, y) = k_1.$$

$$m_1 = \frac{y}{k_1}.$$

And the last condition must be met:

$$\frac{m_3}{k_2} = \frac{m_4}{k_1} = q.$$

The scheme for constructing an Eulerian parallelepiped will look like this:

$$y = (k_1)m_1 \qquad\qquad x = (k_2)m_2$$

$$\downarrow \qquad\qquad\qquad \downarrow \qquad\qquad (4)$$

$$m_3 = k_2 q \qquad\qquad m_4 = k_1 q$$

$$z = k_1 k_2 q$$

Thus, for the selected pair of legs $(x, y)$ from the primitive Pythagorean triple $(x, y, a)$, we construct a system of equations:

$$\begin{cases} y^2 + z^2 = b^2 \\ x^2 + z^2 = c^2 \end{cases}$$

Based on the construction scheme described above, we can build an alternative Eulerian parallelepiped for each obtained Eulerian parallelepiped if we multiply legs $x$ and $y$ by the value $q$. Then the construction scheme will change to the following:

$$qy = (m_1)k_1 q \qquad\qquad qx = (m_2)k_2 q$$

$$\downarrow \qquad\qquad\qquad \downarrow \qquad (5)$$

$$m_2 \qquad\qquad\qquad m_1$$

$$z = m_1 m_2$$

The new coefficients of truncation are written in brackets here. The products $k_1 q$ and $k_2 q$ become new truncated legs, which we will represent as gnomons, placed on the corresponding squares $m_2$ and $m_1$ inside the square lattice with the side $z = m_1 m_2$. Thus, an alternative algorithm for constructing an Eulerian parallelepiped consists of replacing the square of a leg with a gnomon and, conversely, a gnomon with a square of a leg in primitive Pythagorean triples $(m_1, m_3)$ and $(m_2, m_4)$. As a result, we have the following system of equations:

$$\begin{cases} (xq)^2 + (yq)^2 = (aq)^2 \\ (yq)^2 + (m_1 m_2)^2 = (b^*)^2 \\ (xq)^2 + (m_1 m_2)^2 = (c^*)^2 \end{cases}$$

We have described an algorithm for obtaining an Eulerian parallelepiped for a variant in which both pairs of legs $(y, z)$ and $(x, z)$ are constructed as general Pythagorean triples; that is, each pair legs has its own common multiplier. For each Eulerian parallelepiped obtained with this variant, we have obtained a scheme for constructing an alternative Eulerian parallelepiped to it. An example of construction is described in detail in Appendix 2.

Consider a variant of constructing an Eulerian parallelepiped in which the pair of legs $(y, z)$ is constructed as a general Pythagorean triple with a common factor $k_1$, and the second pair of legs $(x, z)$ is constructed as a primitive Pythagorean triple. In this case, $k_2 = 1$.

In this case, we will have the following scheme for finding a common leg $z$:

$$y = (k_1)m_1 \qquad\qquad x = m_2$$

$$\downarrow \qquad\qquad \downarrow$$

$$m_3 = q \qquad\qquad m_4 = k_1 q$$

$$z = k_1 q$$

At the same time, we obtain a scheme for constructing $z$ (formula 4).

For each Eulerian parallelepiped obtained with this variant, we have obtained a scheme for constructing an alternative Eulerian parallelepiped to it. If we multiply legs $x$ and $y$ by the value of $q$, then the $z$ construction scheme will have the form:

$$qy = (m_1)k_1 q \qquad\qquad qx = (m_2)q$$

$$\downarrow \qquad\qquad \downarrow$$

$$m_2 = x \qquad\qquad m_1$$

$$z = m_1 m_2 = m_1 x$$

In this case, we obtain a construction scheme (formula 5).

In addition to constructing Eulerian parallelepipeds from the table of primitive Pythagorean triples and alternative Eulerian parallelepipeds with $x$ and $y$ multiplied by $q$, we can multiply the legs $x$ and $y$ by other coefficients different from $q$. Consider a variant of such a construction of an Eulerian parallelepiped. It occurs if a pair of legs $x$ and $y$ contains multipliers of the following form:

$$y = (r+1)(2r-1) \qquad x = r(2r+3)$$

Let $k_1 = r+1$ and $k_2 = r$. Then we can take the arithmetic mean of the multipliers $(2r-1)$ and $(2r+3)$. This will be number $u = 2r+1$.

We multiply $y$ and $x$ by $(2r+1)$, find the parameter values, and check the conditions necessary and sufficient to construct an Eulerian parallelepiped:

$$k_1 = (r+1); \quad m_1 = (2r-1)(2r+1); \quad l_1 = 2r-1; \quad t_1 = 1.$$

$$m_3 = 2 \cdot 2r = 4r.$$

$$k_2 = r; \quad m_2 = (2r+1)(2r+3); \quad l_2 = 2r+1; \quad t_2 = 1.$$

$$m_4 = 2 \cdot (2r+2) = 4r+4 = 4(r+1).$$

$$\text{GCD}(m_3, x) = k_2 = r;$$

$$\text{GCD}(m_4, y) = k_1 = r+1.$$

$$\frac{m_3}{k_2} = \frac{m_4}{k_1} = q.$$

All three conditions are met and $z = k_1 k_2 q = (r+1)rq$.

Thus, we have constructed an Eulerian parallelepiped of the form:

$$\big((2r+1)y, (2r+1)x, z\big).$$

The common multiplier for the legs $x$ and $y$ is the number $(2r+1)$.

The construction scheme will be as follows:

$$(2r + 1)y = (r + 1)(2r - 1)(2r + 1); \quad (2r + 1)x = (r)(2r + 1)(2r + 3)$$

$$\downarrow \qquad\qquad\qquad\qquad \downarrow$$

$$m_3 = 4r \qquad\qquad\qquad m_4 = 4(r + 1)$$

$$z = 4r(r + 1).$$

Here $m_3$ and $m_4$ are the sides of the inner squares inside the square lattice with the side $z$. Legs squares $m_1 = (2r - 1)(2r + 1)$ and $m_2 = (2r + 1)(2r + 3)$, in the form of corresponding gnomons, are placed on these inner squares and then are assembled into gnomons using the lattice:

$$G_y^* = (2r + 1)y \quad \text{и} \quad G_x^* = (2r + 1)x.$$

Alternative Eulerian parallelepipeds are constructed to the scheme (formula 5).

Reducing by a common multiplier we obtain an alternative Eulerian parallelepiped in the form:

$$4y, 4x, z = (2r - 1)(2r + 1)(2r + 3).$$

We obtain $z = m_1 m_2 = (2r - 1)(2r + 1)^2(2r + 3)$. At the same time, all three legs contain a common multiplier:

$$u = 2r + 1.$$

As a result, any constructed Eulerian parallelepipeds have the same form: a square with a side $z$ and two gnomons are placed on it, which are given by the legs of either a primitive Pythagorean triple or a general Pythagorean triple (scaled).

Herewith, the larger gnomon contains a smaller gnomon inside itself (Fig. 10). The gnomons match with the first terms of the arithmetic progressions describing them.

**Perfect cuboid**

Rectangular parallelepipeds, in which all edges and diagonals of the side faces are expressed in natural numbers, are called Eulerian parallelepipeds.

An Eulerian parallelepiped whose main diagonal is also a natural number (i.e. a parallelepiped, edges, main diagonal and all diagonals of its side faces are natural numbers) is called a perfect cuboid.

In this case, an equation is added to the system of equations (formula 1) which describes the main diagonal (formula 2).

We will consider this equation for the main diagonal; due to additive association, this equation can be written in three different ways:

$$x^2 + (y^2 + z^2) = x^2 + b^2 = d^2;$$

$$y^2 + (x^2 + z^2) = y^2 + c^2 = d^2;$$

$$(x^2 + y^2) + z^2 = a^2 + z^2 = d^2.$$

In general, any Eulerian parallelepiped is represented by a square and two gnomons placed on it. The smaller gnomon is absorbed by the larger one (Fig. 10).

Without loss of generality we will assume that the gnomon $G_y^*$ is smaller than gnomon $G_x^*$.

To solve the equation (formula 2) we need to transform the smaller gnomon and place it to the left of the constructed gnomon $G_x^*$. Let assume that we did it. Denote the newly transformed gnomon as $G_y^{**}$. Consider the constructed square with side $d$ (Fig. 11). Note that in this case we obtained a new transformed gnomon $G_x^{**}$ (Fig. 12).

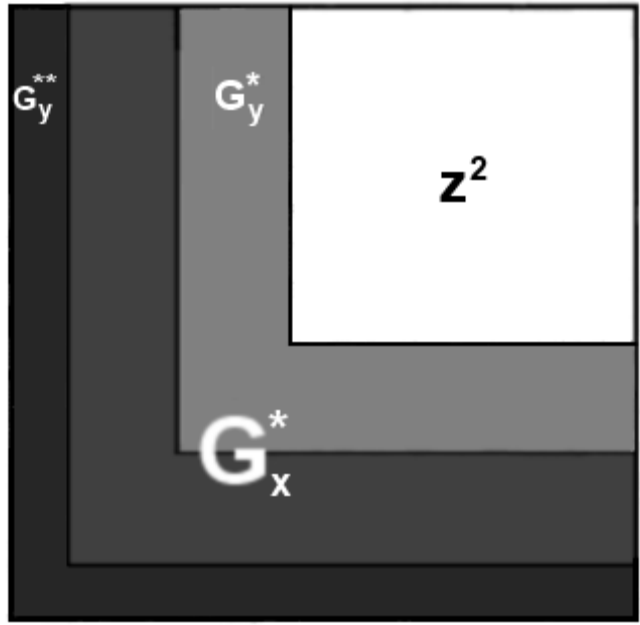
Figure 11

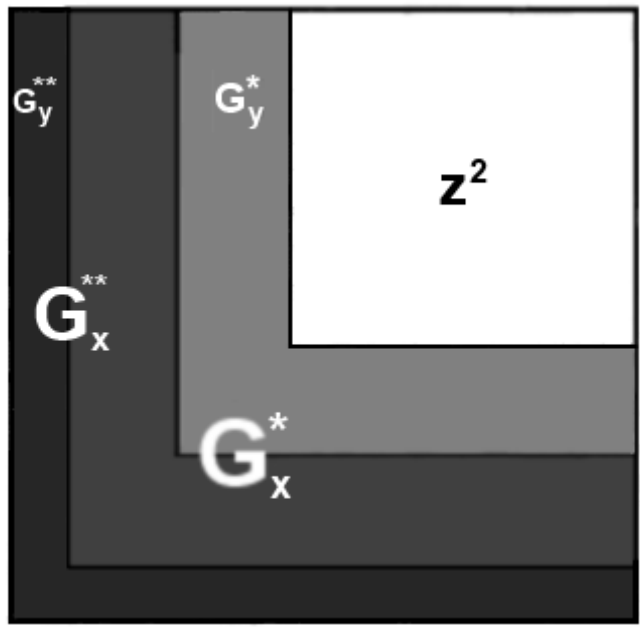
Figure 12

The newly built gnomon $G_y^{**}$ is placing on the square with side $c$ (Fig. 13).

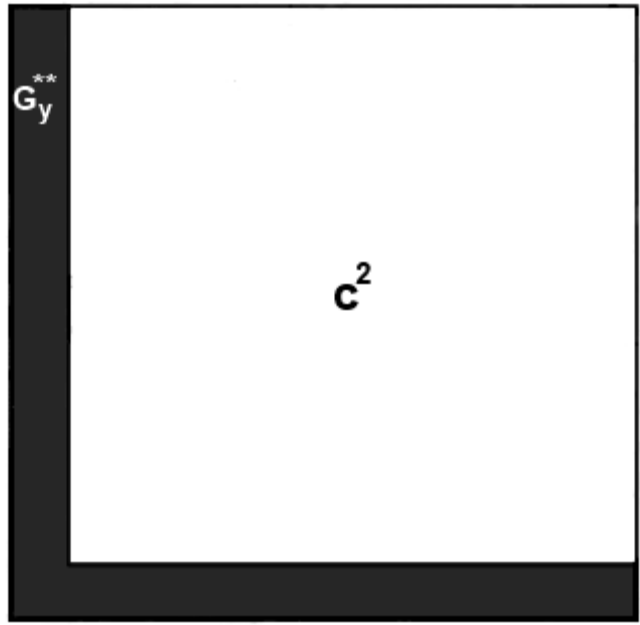

Figure 13

The new transformed gnomon $G_x^{**}$ is placing on the square with side $b$ (Fig. 14).

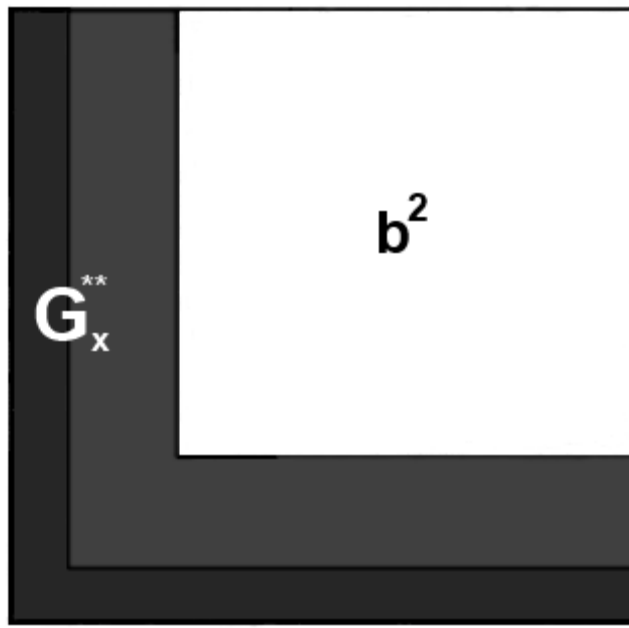

Figure 14

But we obtained two connected gnomons $G_x^{**}$ and $G_y^{**}$ on the left for the squares $x$ and $y$ (Fig. 15).

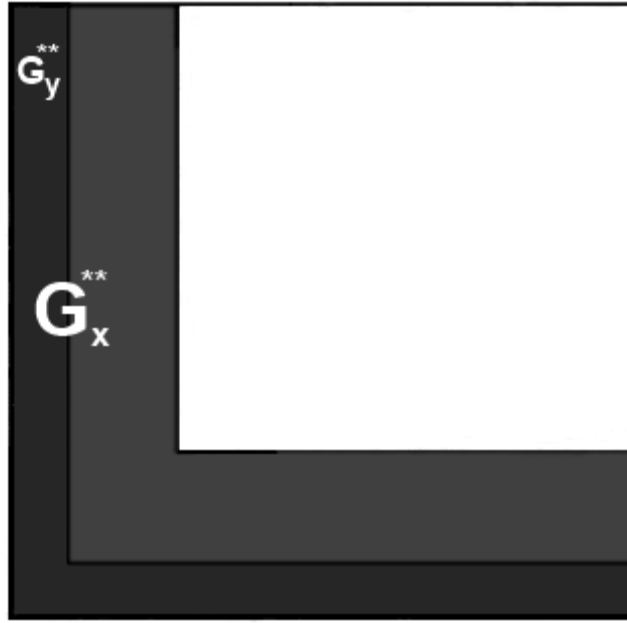

Figure 15

As we described earlier, each of the connected gnomons must place on the square of the other, and have a common outer side equals to $a$ for $x$ and $y$. However, during the construction, the gnomon $G_x^{**}$ is placing on the square with the side $b$ (Fig. 14) and the gnomon $G_y^{**}$ is placing on the square with the side $c$ (Fig.13). In this case:

$$b^2 = y^2 + z^2 > y^2;$$

$$c^2 = x^2 + z^2 > x^2.$$

The desired main diagonal by construction must be equal to $d > a$:

$$d^2 = a^2 + z^2 > a^2.$$

Our construction leads to a contradiction. Consequently, such a transformation of the gnomon $G_y^*$ into the gnomon $G_y^{**}$ cannot exist; therefore, it is impossible to construct an integer main diagonal. This proves the impossibility of constructing a perfect cuboid with integer values of all edges and diagonals.

**Conclusion**

An algorithm for constructing Eulerian parallelepipeds using primitive Pythagorean triples is shown in the paper. The construction of the primitive Pythagorean triples is based on the concept of a generating square with an even side. Using this concept, the relationship of the parameters necessary for the construction of primitive Pythagorean triples is found. Abstract formulas for constructing primitive Pythagorean triples via the traditional, unrelated parameters $m$ and $n$ are replaced by mutually related parameters via the side of the generating square $S = 2tl$, where $t = n; m = l + t$. A sequential increase in the side of the generating square with a constant step equal to 2, starting from 2, leads to the construction of an order on the set of primitive Pythagorean triples using these two related parameters: $S$ and $t$. There is given a table of primitive Pythagorean triples constructed in ascending order of the parameter $S$ and the second-level parameter inside the block $S$, namely $t(S)$.

The paper shows a description of connected gnomons by arithmetic progressions. Using the parameters of the arithmetic progression, it is shown that the connected gnomons place on each other's squares and have a common outer side equal to the hypotenuse.

It is shown that during the transformation of any gnomon (coordinated change in the thickness of the gnomon and the middle term of the arithmetic progression describing the gnomon, provided that the area of the gnomon is preserved), it forms other primitive Pythagorean triples.

Three ways of representing a primitive Pythagorean triple are described: the square of the first leg plus the gnomon of the second; conversely, the square of the second leg plus the gnomon of the first leg; and two connected gnomons. A method of assembling gnomons on a square lattice for legs with a common multiplier is proposed.

The algorithm for constructing Eulerian parallelepipeds is based on the transformation of connected gnomons which leads to the disintegration of an initial taken primitive Pythagorean triple and selection of a new second leg for each transformed gnomon. If the obtained new leg ends up being the same for both transformed gnomons, then this is the solution. The construction is described by a square, the side of which is equal to the new leg, and two transformed gnomons placed on it, with the larger gnomon absorbing the smaller gnomon. At the same time, both transformed gnomons match with the first terms of the arithmetic progressions describing them. In this way they differ from the connected gnomons depicting the primitive Pythagorean triple, which match with the last terms of the arithmetic progressions describing them. Necessary and sufficient conditions for constructing an Eulerian parallelepiped are found.

The algorithm is based on sorting primitive Pythagorean triples from the table of primitive Pythagorean triples and constructing the desired third leg by transforming the original gnomons built on the squares of the legs of a primitive Pythagorean triple with the selection of transformation parameters that meet the necessary and sufficient conditions for constructing an Eulerian parallelepiped. Additionally, an algorithm for constructing an alternative Eulerian parallelepiped for an already constructed basic Eulerian parallelepiped is shown. The construction of an Eulerian parallelepiped for the legs of a primitive Pythagorean triple, which have a special form that allows them to be multiplied by a common multiplier and build an Eulerian parallelepiped by calculating this multiplier, is shown. As a result of the described possible constructions of the Eulerian parallelepiped, we have the same type of construction of the third leg: this is the square of a new leg and two gnomons placed on it, described by the given legs of a primitive Pythagorean triple or by the same legs multiplied by a common multiplier. At the same time, the larger gnomon always absorbs the smaller gnomon.

Based on the result of constructing an Eulerian parallelepiped we made an attempt to build a perfect cuboid. To accomplish this, it was necessary to solve the equation describing the main diagonal. Assuming that this solution is possible, we transformed the smaller gnomon in such a way that it was transferred to the end of the large gnomon (Fig.11). As a result of this construction, we received two gnomons transformed again; however, via the construction process, the gnomons became the connected gnomons that led to a contradiction in solving the equation for the main diagonal, which indicated the impossibility of constructing a perfect cuboid.

Therefore it is impossible to build a perfect cuboid with integer parameters.

Only two productive constructions using gnomons are possible: the initial position - when the connected gnomons of the primitive Pythagorean triple are overlapped with the last terms of the arithmetic progressions describing them; and the final position, when the transformed gnomons build an Eulerian parallelepiped, while they are overlapping with the first terms of the arithmetic progressions describing them. Other transformations of gnomons do not lead to significant results.

# Appendix 1

# Fragment of a table of primitive Pythagorean triples constructed with increasing parameters $S, t(S)$

| $N.n_i$ | $S$ | $t$ | $l$ | $x = S + l^2$ | $y = S + 2t^2$ | $a = S + l^2 + 2t^2$ |
|---|---|---|---|---|---|---|
| 1.1 | 2 | 1 | 1 | 3 | 4 | 5 |
| 2.1 | 4 | 2 | 1 | 5 | 12 | 13 |
| 3.1 | 6 | 1 | 3 | 15 | 8 | 17 |
| 3.2 |  | 3 | 1 | 7 | 24 | 25 |
| 4.1 | 8 | 4 | 1 | 9 | 40 | 41 |
| 5.1 | 10 | 1 | 5 | 35 | 12 | 37 |
| 5.2 |  | 5 | 1 | 11 | 60 | 61 |
| 6.1 | 12 | 2 | 3 | 21 | 20 | 29 |
| 6.2 |  | 6 | 1 | 13 | 84 | 85 |
| 7.1 | 14 | 1 | 7 | 63 | 16 | 65 |
| 7.2 |  | 7 | 1 | 15 | 112 | 113 |
| 8.1 | 16 | 8 | 1 | 17 | 144 | 145 |
| 9.1 | 18 | 1 | 9 | 99 | 20 | 101 |
| 9.2 |  | 9 | 1 | 19 | 180 | 181 |
| 10.1 | 20 | 2 | 5 | 45 | 28 | 53 |
| 10.2 |  | 10 | 1 | 21 | 220 | 221 |
| 11.1 | 22 | 1 | 11 | 143 | 24 | 145 |
| 11.2 |  | 11 | 1 | 23 | 264 | 265 |
| 12.1 | 24 | 4 | 3 | 33 | 56 | 65 |
| 12.2 |  | 12 | 1 | 25 | 312 | 313 |
| 13.1 | 26 | 1 | 13 | 195 | 28 | 197 |
| 13.2 |  | 13 | 1 | 27 | 364 | 365 |
| 14.1 | 28 | 2 | 7 | 77 | 36 | 85 |
| 14.2 |  | 14 | 1 | 29 | 420 | 421 |
| 15.1 | 30 | 1 | 15 | 255 | 32 | 257 |
| 15.2 |  | 3 | 5 | 55 | 48 | 73 |
| 15.3 |  | 5 | 3 | 39 | 80 | 89 |
| 15.4 |  | 15 | 1 | 31 | 480 | 481 |
| 16.1 | 32 | 16 | 1 | 33 | 544 | 545 |
| 17.1 | 34 | 1 | 17 | 223 | 36 | 225 |
| 17.2 |  | 17 | 1 | 35 | 612 | 613 |
| 18.1 | 36 | 2 | 9 | 117 | 44 | 125 |
| 18.2 |  | 18 | 1 | 37 | 684 | 685 |
| 19.1 | 38 | 1 | 19 | 399 | 40 | 401 |
| 19.2 |  | 19 | 1 | 39 | 760 | 761 |
| 20.1 | 40 | 4 | 5 | 65 | 72 | 97 |
| 20.2 |  | 20 | 1 | 41 | 840 | 841 |
| 21.1 | 42 | 1 | 21 | 483 | 44 | 485 |
| 21.2 |  | 3 | 7 | 91 | 60 | 109 |
| 21.3 |  | 7 | 3 | 51 | 140 | 149 |
| 21.4 |  | 21 | 1 | 43 | 924 | 925 |
| 22.1 | 44 | 2 | 11 | 165 | 52 | 173 |
| 22.2 |  | 22 | 1 | 45 | 1012 | 1013 |
| 23.1 | 46 | 1 | 23 | 575 | 48 | 577 |
| 23.2 |  | 23 | 1 | 47 | 1104 | 1105 |
| 24.1 | 48 | 8 | 3 | 57 | 176 | 185 |
| 24.2 |  | 24 | 1 | 49 | 1200 | 1201 |
| 25.1 | 50 | 1 | 25 | 675 | 52 | 677 |
| 25.2 |  | 25 | 1 | 51 | 1300 | 1301 |

| $N.n_i$ | $S$ | $t$ | $l$ | $x = S + l^2$ | $y = S + 2t^2$ | $a = S + l^2 + 2t^2$ |
|---|---|---|---|---|---|---|
| 26.1 | 52 | 2 | 13 | 221 | 60 | 229 |
| 26.2 |  | 26 | 1 | 53 | 1404 | 1405 |
| 27.1 | 54 | 1 | 27 | 783 | 56 | 785 |
| 27.2 |  | 27 | 1 | 55 | 1512 | 1513 |
| 28.1 | 56 | 4 | 7 | 105 | 88 | 137 |
| 28.2 |  | 28 | 1 | 57 | 1624 | 1625 |
| 29.1 | 58 | 1 | 29 | 899 | 60 | 901 |
| 29.2 |  | 29 | 1 | 59 | 1740 | 1741 |
| 30.1 | 60 | 2 | 15 | 285 | 68 | 293 |
| 30.2 |  | 6 | 5 | 85 | 132 | 157 |
| 30.3 |  | 10 | 3 | 69 | 260 | 269 |
| 30.4 |  | 30 | 1 | 61 | 1860 | 1861 |
| 31.1 | 62 | 1 | 31 | 1023 | 64 | 1025 |
| 31.2 |  | 31 | 1 | 63 | 1984 | 1985 |
| 32.1 | 64 | 32 | 1 | 65 | 2112 | 2113 |
| 33.1 | 66 | 1 | 33 | 1155 | 68 | 1157 |
| 33.2 |  | 3 | 11 | 187 | 84 | 205 |
| 33.3 |  | 11 | 3 | 75 | 308 | 317 |
| 33.4 |  | 33 | 1 | 67 | 2244 | 2245 |
| 34.1 | 68 | 2 | 17 | 357 | 76 | 365 |
| 34.2 |  | 34 | 1 | 69 | 2380 | 2381 |
| 35.1 | 70 | 1 | 35 | 1295 | 72 | 1297 |
| 35.2 |  | 5 | 7 | 119 | 120 | 169 |
| 35.3 |  | 7 | 5 | 95 | 168 | 193 |
| 35.4 |  | 35 | 1 | 71 | 2520 | 2521 |
| 36.1 | 72 | 4 | 9 | 153 | 104 | 185 |
| 36.2 |  | 36 | 1 | 73 | 2664 | 2665 |
| 37.1 | 74 | 1 | 37 | 1443 | 76 | 1445 |
| 37.2 |  | 37 | 1 | 75 | 2812 | 2813 |
| 38.1 | 76 | 2 | 19 | 437 | 84 | 445 |
| 38.2 |  | 38 | 1 | 77 | 2964 | 2965 |
| 39.1 | 78 | 1 | 39 | 1599 | 80 | 1601 |
| 39.2 |  | 3 | 13 | 247 | 96 | 265 |
| 39.3 |  | 13 | 3 | 87 | 416 | 425 |
| 39.4 |  | 39 | 1 | 79 | 3120 | 3121 |
| 40.1 | 80 | 8 | 5 | 125 | 208 | 233 |
| 40.2 |  | 40 | 1 | 81 | 3280 | 3281 |
| 41.1 | 82 | 1 | 41 | 1763 | 84 | 1765 |
| 41.2 |  | 41 | 1 | 83 | 3444 | 3445 |
| 42.1 | 84 | 2 | 21 | 525 | 92 | 533 |
| 42.2 |  | 6 | 7 | 133 | 156 | 205 |
| 42.3 |  | 14 | 3 | 93 | 476 | 485 |
| 42.4 |  | 42 | 1 | 85 | 3612 | 3613 |
| 43.1 | 86 | 1 | 43 | 1935 | 88 | 1937 |
| 43.2 |  | 43 | 1 | 87 | 3784 | 3785 |
| 44.1 | 88 | 4 | 11 | 209 | 120 | 241 |
| 44.2 |  | 44 | 1 | 89 | 3960 | 3961 |
| 45.1 | 90 | 1 | 45 | 2115 | 92 | 2117 |
| 45.2 |  | 5 | 9 | 171 | 140 | 221 |
| 45.3 |  | 9 | 5 | 115 | 252 | 277 |
| 45.4 |  | 45 | 1 | 91 | 4140 | 4141 |
| 46.1 | 92 | 2 | 23 | 621 | 100 | 629 |
| 46.2 |  | 46 | 1 | 93 | 4324 | 4325 |
| 47.1 | 94 | 1 | 47 | 2303 | 96 | 2305 |
| 47.2 |  | 47 | 1 | 95 | 4512 | 4513 |
| 48.1 | 96 | 16 | 3 | 105 | 608 | 617 |

| $N.n_i$ | $S$ | $t$ | $l$ | $x = S + l^2$ | $y = S + 2t^2$ | $a = S + l^2 + 2t^2$ |
|---|---|---|---|---|---|---|
| 48.2 |  | 48 | 1 | 97 | 4704 | 4705 |
| 49.1 | 98 | 1 | 49 | 2499 | 100 | 2501 |
| 49.2 |  | 49 | 1 | 99 | 4900 | 4901 |
| 50.1 | 100 | 2 | 25 | 725 | 108 | 733 |
| 50.2 |  | 50 | 1 | 101 | 5100 | 5101 |

**Appendix 2**

**Construction of a minimal Eulerian parallelepiped**

Consider the dynamic discrete transformation of gnomons. Take the primitive Pythagorean triple (44, 117, 125). The even leg is 44. The odd leg is equal to 117.

The initial position of the gnomon corresponding to the square of an even leg is as follows:

$$2t(l + t);$$

$$44 = 2 \cdot 2 \cdot 11;$$

$$t = 2; l = 11 - 2 = 9; S = 2 \cdot 2 \cdot 9 = 36;$$

The thickness of the gnomon is $T_y = 2t^2$; $y = S + T_y = 36 + 2 \cdot 2^2 = 44$. The thickness of the gnomon is equal to the number of terms of the arithmetic progression. The outer side of the gnomon is equal to the hypotenuse:

$$a = y + l^2 = 44 + 81 = 125.$$

The gnomon is placing on a square with side x=117. The side of an odd square can be calculated by the formula:

$$\frac{2(t+l)^2 - 2t^2}{2} = (t+l)^2 - t^2 = (t+l-t) \cdot (t+l+t) = l \cdot (l+2t).$$

The middle term of the arithmetic progression for the gnomon corresponding to an even square (a square with an even side) is $2 \cdot (t+l)^2$. Here the obtained formula $l \cdot (l + 2t)$ corresponds to the formula of an odd leg.

$$l \cdot (l + 2t) = 9 \cdot (9 + 2 \cdot 2) = 9 \cdot 13 = 117.$$

The initial position of the gnomon corresponding to the square of an odd leg is as follows:

$$l \cdot (l + 2t);$$

The parameters $S, l, t$ are common to all numbers of this primitive Pythagorean triple.

The thickness of the gnomon is $T_x = l^2; x = S + T_x = 36 + 9^2 = 117$. The thickness of the gnomon is equal to the number of terms of the arithmetic progression describing this gnomon. The outer side of the gnomon is equal to the hypotenuse:

$$a = x + 2t^2 = 117 + 8 = 125.$$

The gnomon is placing on a square with side $y = 44$. The side of an even square can be calculated by the formula:

$$\frac{(l + 2t)^2 - l^2}{2} = \frac{(l + 2t - l) \cdot (l + 2t + l)}{2} = 2t \cdot (l + t).$$

The resulting formula $2t \cdot (l + t)$ corresponds to the formula of an even leg.

The middle term of the arithmetic progression for the gnomon corresponding to the square of an odd number is equal to $(l + 2t)^2$.

Here we have shown that each gnomon placed on the square of its paired leg.

Let's write down the arithmetic progressions corresponding to both gnomons.

The first term of the arithmetic progression for the gnomon corresponding to the square of $x$ is equal to

$$2 \cdot 44 + 1 = 89.$$

The number of terms is $l^2 = 9^2 = 81$.

Let's write down all the terms of this progression:

89 91 93 95 97 99 101 103 105 107

109 111 113 115 117 119 121 123 125 127

129 131 133 135 137 139 141 143 145 147

149 151 153 155 157 159 161 163 165 167

169 171 173 175 177 179 181 183 185 187

189 191 193 195 197 199 201 203 205 207

209 211 213 215 217 219 221 223 225 227

229 231 233 235 237 239 241 243 245 247

249

The middle term of the arithmetic progression is equal to $(l+2t)^2 = 13^2 = 169$. We highlighted it in color.

The first term of the arithmetic progression for the gnomon corresponding to the square of $y$ is equal to

$$2 \cdot 117 + 1 = 235.$$

The number of terms is $2t^2 = 2 \cdot 2^2 = 8$.

Let's write down all the terms of this progression:

235 237 239 241 243 245 247 249

The middle term of the arithmetic progression is equal to $2 \cdot (l+t)^2 = 2 \cdot 11^2 = 242$

In this arithmetic progression the middle term and the number of terms are even numbers.

All terms of the arithmetic progression with a smaller amount of terms match with the last terms in the arithmetic progression with a larger amount of terms. Here we highlighted it in green color.

General Pythagorean triples with a common coefficient $k$ have an increase in the number of terms and middle terms of the arithmetic progression by $k$ times.

To build an Eulerian parallelepiped the gnomons of both legs are transformed in such a way that they can be placed on the same square.

We transform an even leg into the product of two multipliers

$$y = k_1 m_1.$$

The truncation coefficient of an even leg has the form $k_1 = 4i$, $i \in \mathbb{N}$. The second cofactor is $m_1 \neq 1$. And if it is even, it has a factor equal to four.

In our case $y = 44 = 4 \cdot 11$. $k_1 = 4$; $m_1 = 11$.

The second cofactor is an odd number. So this is an odd leg, which is represented by the formula $l_1 \cdot (l_1 + 2t_1)$. Since this cofactor is a prime number, then $l_1 = 1$. It follows $t_1 = 5$.

The paired leg for $m_1$, which is included with it in the primitive Pythagorean triple, is calculated by the formula $m_3 = 2t_1 \cdot (l_1 + t_1)$. It is equal to $2 \cdot 5 \cdot (1 + 5) = 60$.

In the next step, we find the truncation coefficient for the odd leg $x$:

$$k_2 = \text{GCD}(m_3, x).$$

$$m_3 = 60 = 3 \cdot 4 \cdot 5; \quad x = 3 \cdot 3 \cdot 13; \quad k_2 = 3.$$

Represent the leg $x$ as a product of two cofactors $x = k_2 \cdot m_2$.

In turn, $m_2 = 3 \cdot 13$. Since this leg is odd, it is represented by the formula $m_2 = l_2 \cdot (l_2 + 2t_2)$. Select $l_2$ from the list of possible representatives $\{1, 3\}$. For $l_2 = 3$ we have $t_2 = 5$.

The paired leg for $m_2$, which is included with it in the primitive Pythagorean triple, is calculated by the formula $m_4 = 2t_2 \cdot (l_2 + t_2)$. It is equal to $2 \cdot 5 \cdot (3 + 5) = 80$.

We have constructed the legs $m_3 = k_2 \cdot 20$; $m_4 = k_1 \cdot 20$. Let's denote the common multiplier of these legs as $q$. This will be the side of the smallest square of the transformational square lattice. The full side of the lattice will be equal to $z = k_1 k_2 q$. The side of the transformational square lattice will be the desired third leg of the Eulerian parallelepiped:

$$\begin{cases} y^2 + x^2 = a^2 \\ y^2 + z^2 = b^2 \\ x^2 + z^2 = c^2 \end{cases}$$

We first select a square with the side $m_3 = k_2 q$ in the square lattice. The number of such squares will be equal to $k_1^2$. We present in the form of a gnomon $m_1^2$. Such a gnomon is placing on every square $m_3$ in a square lattice. Then all the gnomons are combined outside the lattice t. i. they are assembled and placed on the left on the square lattice. In this case, the thickness of the total gnomon $y$ becomes equal to $T_y^* = k_1 l_1^2$. Thus, the gnomon was transformed to fit on the square with the $z$ side. The outer side of the gnomon is equal to the hypotenuse

$$b = z + k_1 l_1^2$$

We next select a square with the side $m_4 = k_1 q$. The number of such squares will be equal to $k_2^2$. We present $m_2^2$ in the form of a gnomon. Such a gnomon is placing on every square $m_4$ in a square lattice. Then all the gnomons are combined outside the lattice, i. e. they are assembled and placed on the left on the square lattice. In this case, the thickness of the total gnomon $x$ becomes equal to $T_x^* = k_2 l_2^2$. Thus, the gnomon was transformed to fit on the square with the $z$ side. The outer side of the gnomon is equal to the hypotenuse

$$c = z + k_2 l_2^2$$

Similarly, we can represent this as a transformation of arithmetic progressions for gnomons $x$, $y$. This is constructed as the difference between the middle terms of arithmetic progressions and the number of their terms when finding a common first term. Since the placing of gnomons occurs on the same square.

Let's find the middle terms of both arithmetic progressions:

$$s_y = \frac{y^2}{T_y^*} = \frac{y^2}{k_1 l_1^2} = \frac{4^2 \cdot 11^2}{4 \cdot 1^2} = 4 \cdot 11^2 = 484;$$

$$s_x = \frac{x^2}{T_x^*} = \frac{x^2}{k_2 l_2^2} = \frac{3^2 \cdot 3^2 \cdot 13^2}{3 \cdot 3^2} = 3 \cdot 13^2 = 507.$$

Let's find the side of the square on which both gnomons are placing:

$$z = \frac{s_y - T_y^*}{2} = \frac{s_x - T_x^*}{2};$$

$$z = \frac{484 - 4}{2} = \frac{507 - 27}{2} = 240.$$

We obtained the same square on which both gnomons were placing. Thus, by transforming gnomons and, accordingly, describing their arithmetic progressions, we solved the problem of finding the third leg $z$ necessary for constructing an Eulerian parallelepiped.

Let's write down the arithmetic progressions corresponding to both transformed gnomons.

Let's write down the arithmetic progressions corresponding to both transformed gnomons.

The first term of the arithmetic progression, the same for both gnomons, is equal to $2z + 1$. The number of terms of the arithmetic progression corresponding to the square of $y$ is equal to $T_y^* = 4$.

Let's write down all the terms of this progression:

<span style="color:green">481 483 485 487</span>

The middle term of the arithmetic progression is 484.

The number of terms of the arithmetic progression corresponding to the square of $x$ is $T_x^* = 27$.

Let's write down all terms of this progression:

<span style="color:green">481 483 485 487</span> 489 491 493 495 497 499

501 503 505 <span style="color:red">507</span> 509 511 513 515 517 519

521 523 525 527 529 531 533

The middle term is equal to 507. It is highlighted in red. The arithmetic progression with smaller amount of terms is the initial part of the arithmetic progression with larger amount of terms. Common terms are highlighted in green.

Thus, we have constructed an Eulerian parallelepiped:

$$(44, 117, 240).$$